\documentclass[12pt,a4paper]{article}
\usepackage[latin1]{inputenc}
\usepackage{amsmath}
\usepackage{amssymb}
\usepackage{amsthm}
\usepackage{graphicx}
\usepackage{color}
\usepackage{fancybox}
\usepackage{amsfonts}
\usepackage{pb-diagram}
\usepackage{float}

\newtheorem{theo}{Theorem}
\newtheorem{prop}[theo]{Proposition}
\newtheorem{defi}{Definition}
\newtheorem*{re}{Remark}

\newtheorem{lemme}[theo]{Lemma}

\theoremstyle{definition}
\newtheorem*{dem}{Proof}

\renewcommand{\textbf}[1]{\begingroup\bfseries\mathversion{bold}#1\endgroup}

\newcommand{\ffp}{\mbox{forest-fire process}}
\newcommand{\ffps}{\mbox{forest-fire processes}}

\def\zd{\mathbb{Z}^{d}}
\def\bx{\mathcal{B}(X)}

\usepackage{tikz}
\usetikzlibrary{calc,trees,positioning,arrows,chains,shapes.geometric,%
    decorations.pathreplacing,decorations.pathmorphing,shapes,%
    matrix,shapes.symbols}

\tikzset{   
terminal/.style={
rectangle,minimum size=6mm,rounded corners=3mm,
very thick,draw=black!50,
top color=white,bottom color=black!20,
},
nonterminal/.style={
rectangle,
minimum size=6mm,
very thick,
draw=blue!50!black!50, 
top color=white, 
bottom color=blue!50!black!20, 
font=\itshape
},
}  

\begin{document}

\begin{center}
{\large Existence of a stationary distribution for multi-dimensional infinite volume forest-fire processes} 

\vspace*{0.2cm}

A. Stahl

\vspace*{0.2cm}

University of Toulouse, France
\end{center}
\vspace*{1cm}

\textbf{\textit{Abstract.}} Consider the following forest-fire process on a connected graph. Each site of the graph can be either occupied or vacant. A vacant site becomes occupied with rate $1$. A site is ignited with rate $\lambda$, and its whole occupied cluster burns instantaneously. 

The purpose of this paper is to show the existence of a stationary distribution for forest-fire processes on $\zd$, for $d \geq 2$. We define a distribution on $\left\lbrace 0,1 \right\rbrace ^{\zd} $ as a limit of a sequence of invariant distributions of finite volume forest-fire processes, and then show it is a stationary distribution for $\ffps$ on $\zd$, $d \geq 2$.

\section*{Introduction} 
\label{intro}

Forest-fire models have been introduced within the study of self-organised criticality (\cite{Bak90}). The idea behind self-organised criticality is the evolution of the system toward a critical behavior, only driven by its own interaction rules (see \cite{Jen} for example). In the nineties, a forest-fire model in discrete time and on a finite graph, introduced by Drossel and Schwabl \cite{DroSch92}, has been studied in the physics literature.  In $2005$ in \cite{BeJa} van den Berg and Jarai studied the asymptotic behavior of a continuous-time version of the Drossel-Schwabl forest-fire model, on an infinite graph : $\mathbb{Z}$. Then this model has been studied on $\mathbb{Z}$, on $\zd$ and on connected graphs with a vertex degree bounded by $d$, for $d \geq 2$. 

Let $G=(V,E)$ be a connected graph with vertex degree bounded by $d \geq 2$. $V$ is the set of vertices and $E$ the set of edges of the graph. In the model, each site of $V$ can be vacant or occupied by a tree. Two sites of $V$ are neighbors if they are linked by an edge in $E$. The evolution is driven by the initial configuration and two families of Poisson processes, one governing the birth of trees and the other the ignition of a site. Throughout this paper ``birth" of a tree and ``growth" of a tree have the same meaning. An empty site becomes occupied if and only if there is a birth attempt. An occupied site becomes empty either if the site is ignited or if a site of its connected cluster of occupied sites is ignited. The second case models the propagation of fire. As the burning is instantaneous here, the speed of the fire propagation is infinite. 

More precisely, the \textit{configuration at time} $t$ of the forest is the process $((\eta_{t,x})_{x \in V})_{t \geq 0}$ where $\eta_{t,x}$ takes the value $1$ if the site $x$ is occupied at time $t$ and $0$ otherwise. If the site $x$ is occupied at time $t$, its \textit{connected cluster} $C_{t,x}$ is the set of all occupied sites connected to $x$ by a path of occupied sites. 

\begin{defi}(open path and connected cluster of a site)
\begin{enumerate}
\item A path between 2 sites $x$ and $y$ in $V$ is a set of sites $ \lbrace z_i \rbrace_{1 \leq i \leq n} $ of $V$ such that $z_1=x$, $z_n=y$ and for all integers $i \in \lbrace 1, \dots , n-1 \rbrace $ there is an edge of $E$ linking $z_i$ and $z_{i+1}$. The path is open at time $t$ if all its sites are occupied at time $t$.
\item The connected cluster $C_{t,x}$ of a site $x$ at time $t$ is the set of all sites of $V$ linked to $x$ with an open path.  
\end{enumerate}
\end{defi}

To each site $x \in V$ are associated two Poisson processes, $(G_{t,x})_{t \geq 0}$ and  $(I_{t,x})_{t \geq 0}$, with intensity $1$ and $\lambda>0$ respectively. The process $G_x$ represents the birth attempts on the site $x$ whereas $I_x$ represents ignition attempts on the site $x$. These two processes are independent and moreover independent of each Poisson process of all the other sites. The following definition gives a rigourous definition of the forest-fire process.

\begin{defi}

Let $S$ be a subset of $V$ and $\lambda$ a non negative real number. A $\ffp$ on $S$ with parameter $\lambda>0$ is a process $(\overline{\eta}_{t})_{t \geq 0}$ where $\overline{\eta}_{t}=(\overline{\eta}_{t,x} )_{x \in
S}=(\eta_{t,x},G_{t,x},I_{t,x})_{x \in S}$, with values in $(\{
0,1\}  \times \mathbb{N} \times \mathbb{N})^{S}$, that has the following properties.
\begin{enumerate}
\item{For all sites $x$ in $S$, the processes $(G_{t,x})_{t \geq 0}$ and $(I_{t,x})_{t \geq 0}$,
are independent Poisson processes with parameter $1$ and $\lambda$ respectively.}
\item {For all sites $x$ in $S$, the process $(\eta_{t,x},G_{t,x},I_{t,x})_{t \geq 0}$ is right continous and left limited.
}

\item{For all $x \in S$ and all $t \geq 0$,
\begin{enumerate}
\item{if there is a growth of a tree at the site $x$ at time $t$, then the site $x$ is occupied at time $t$, $$G_{t^{-},x}< G_{t,x}
\Rightarrow \eta_{t,x}=1,$$}
\item{if the site $x$ becomes occupied at time $t$, then a tree must have grown at the site $x$ at time $t$
$$\eta_{t^{-},x}< \eta_{t,x}
\Rightarrow G_{t^{-},x}< G_{t,x},$$ }
\item{if the site $x$ is ignited at time $t$, then all the sites of its connected cluster get vacant at time $t$,
$$I_{t^{-},x}< I_{t,x} \Rightarrow \forall y \in C_{t^{-},x}, \
\eta_{t,y}=0,$$}
\item{if the site $x$ gets vacant at time $t$, then the connected cluster of the site $x$ has been ignited at time $t$, $$\eta_{t^{-},x}> \eta_{t,x} \Rightarrow
\exists y \in C_{t^{-},x} :I_{t^{-},y}< I_{t,y}.$$}
\end{enumerate}
}
\end{enumerate}

\end{defi}

The notation $\eta(t)$ will be used to denote the configuration at time $t$ : $\lbrace \eta_{t,x} \rbrace_{x \in S}$.

\vspace{0.3cm}

On a finite graph, the forest-fire process is uniquely determined by the following contruction. At time zero, the initial configuration of the forest is given. Then, the evolution is governed by the Poisson processes. Due to the finite number of sites, all the Poissonian events can be ordered in time. Then, four situations can occur. Imagine that there is a growth attempt on the site $x$. If the site $x$ is vacant, then it becomes occupied. But if $x$ is already occupied, nothing happens. Imagine now that there is an ignition attempt on the site $y$. If the site $y$ is vacant, nothing happens. But if $y$ is occupied, then its whole connected cluster burns down immediately. At each Poissonian event, the configuration is updated using these rules. 

\vspace{0.3cm}

However on an infinite graph, infinitely many Poissonian events can occur in a finite interval of time. Thus the construction given above does not work anymore. Does this forest-fire process exist on infinite graphs? The answer is yes. It can easily be seen on $\mathbb{Z}$ using a graphical construction (\cite{Lig}) but it is harder on $\zd$ for $d \geq 2$. On $\mathbb{Z}$, if an initial configuration has infinitely many empty sites, then at a time $t>0$ there are almost surely infinitely many empty sites that have been empty during the whole interval of time $[0,t]$. An empty site blocks the propagation of fire on $\mathbb{Z}$. So, $\mathbb{Z}$ can be divided into finite pieces and the previous construction can be used on each piece. But on $\zd$ for $d \geq 2$, an empty site alone cannot block fires. By using a sequence of finite-volume processes and finding a control of the influence of the state of long-distance sites, Dürre has shown in \cite{Dur-exist} the existence of forest-fire processes on $\zd$ for all initial configurations that contain no infinite clusters. This result can be extended to connected graphs with vertex degree bounded by $d \geq 2$. However, nothing ensures that given an initial configuration and all the Poisson processes, the forest-fire process is unique. This uniquenesss result is proven for all parameter $\lambda>0$ in \cite{Dur-these}, with an assumption on the size of the clusters in the initial configuration. The question of uniqueness for all initial configurations is still open.

In this paper, we will be interested in stationary distributions of the forest-fire process. The $\ffp$ is a Markov process but not a Feller process. Thus the usual arguments used with interacting particles systems will not work here (see for example \cite{Lig}). Brouwer and Pennanen have shown in {\cite{BrPen} the existence of at least one stationary distribution for the process on $\mathbb{Z}$. In \cite{BrFo-uni} Bressaud and Fournier have shown that the stationary distribution is unique when the parameter $\lambda$ is equal to $1$. 

The purpose of this note is to prove the existence of a stationary distribution for $\ffps$ on $\zd$ with $d$ larger than $2$. The method of Brouwer and Pennanen combined with some tools introduced by Dürre in his thesis \cite{Dur-these} will be used. Since the existence of forest-fire processes on $\zd$ is known yet for initial configurations that contain no infinite clusters, we will consider only such forest-fire processes in this paper.

\begin{theo}
\label{main th}
A $\ffp$ on $\zd$ with $d \geq 2$, with parameter $\lambda>0$ has at least one stationary translation-invariant distribution.
\end{theo}

This paper is divided into three parts. The first section presents the main notations and the tools introduced by Dürre in \cite{Dur-these} that will be used in this paper. Then, the construction of the candidate for the stationary distribution is explained in the second part. Finally the last section is devoted to the proof of Theorem~\ref{main th}. 

\section{Notations and tools}

The purpose of this section is to present the main notations and the tools introduced by Dürre in \cite{Dur-these}.

\subsection{Notations} 

\label{notations} 

Unless stated otherwise, in this section, $G$ will denote a connected graph with a vertices set $V$ and an edges set $E$. A bound on the vertex degree of the graph $G$ will be denoted by $d_G$.

\vspace{0.3cm}

For $x=(x_1,\dots,x_d) \in \zd$, $|x|_\infty=\sup \left \lbrace |x_i|, i=1 \dots d\right \rbrace $ and $\displaystyle{|x|_1= \sum_{i=1}^{d} |x_i|} $.

\vspace{0.3cm}

The expression \textit{the graph $\zd$} will refer to a graph with vertices set $V^{(d)}=\zd$ and edges set $E^{(d)}$ corresponding to the hypercubic lattice : $E^{(d)}= \big \lbrace \lbrace x,y \rbrace \in \zd \times \zd : |x-y|_1 =1 \big \rbrace$.

\vspace{0.3cm}

The site $(0,\dots,0) \in \zd$ will be called the origin. We will work with specifc finite subsets of $\zd$ called boxes, and defined by \[ B_k= \{ x \in \zd : |x|_\infty \leq k  \}.\] Here $k$ is called the radius of the box $B_k$. 

\vspace{0.3cm}

The set of neighbors of a finite set of sites $S$ is called the boundary of $S$ : \[N(S)=\lbrace x \in V \setminus S , \exists \ y \in S  \mbox{ linked by an edge in } E \mbox{ to } x \rbrace.\] We can define the set $\overline{S}=S \cup N(S)$.

\vspace{0.3cm}

Let $X= \left\lbrace 0,1 \right\rbrace ^{\zd} $ be the set of all possible configuations and $\bx$ the associated borel $\sigma$-field. Let $\mu$ be a measure on $(X,\bx)$. For a finite subset $J$ of $\zd$ we denote by $\mu_{|J}$ the restriction of $\mu$ to $J$. For a $\ffp$ $\overline{\eta}$ with an initial configuration with law $\mu$, let us denote by $P^{\mu}(\eta_t \in \cdot )$ its distribution at time $t$.

\vspace{0.3cm}

\subsection{Dürre's tools}

\label{tools} 

\vspace{0.3cm}

In this model, the problem is to control the influence of long-distance sites. Two notions introduced by Dürre in his thesis will be used in this paper. The first one is a condition on the tail of the distribution of cluster size, named $\mbox{CCSB}$. The second one is a process, called the blur-process, used to study the influence on a given site of long-distant sites. Throughout this section, $\ffps$ on connected graphs will be considered.

\vspace{0.3cm}

\subsubsection{CCSB condition}

\vspace{0.3cm}

An important quantity to consider here is the size of clusters. In his thesis, Dürre restricted to $\ffps$ with initial configurations satisfying what he called the conditionned cluster size bound condition to prove the uniqueness. In order to use his results, we will need to use this conditionned cluster size bound condition. 

This condition is a constraint on the cluster size distribution. It implies a uniform bound on the cluster size distribution. It is stronger in the sense that the result still holds when conditionning by the event that a specific configuration of a finite number of sites outside the cluster occurs.

\begin{defi}(Conditionned cluster size bound, CCSB)
For all $s \geq 0$, $\delta \geq 0$ and $m \in \mathbb{N}$, $\overline{\eta}$ has $\mbox{CCSB}(s,\delta,m)$ if it has the following property. Let $B$ an $D$ be two finite subsets of $V$ and $x \in V \setminus D$. Then conditionned on the occurence of the event $\displaystyle{\cup_{y \in B} C_{s,y}=D}$, the probability that the size of the cluster at $x$ is larger than $m$ at time $s$ is smaller than or equal to $\delta$. And almost surely the cluster at $x$ is finite at time $s$. 

More formally, for all finite $B,D \subset V$, for all $x \in V \setminus D$, 
\begin{equation}
 P \big(  \big \lbrace |C_{s,x}|>m \big \rbrace  \bigcap \big \lbrace \bigcup_{y \in B} C_{s,y}=D \big \rbrace \big)  \leq \delta \hspace*{0.1cm} P \big(  \bigcup_{y \in B} C_{s,y}=D \big) 
\end{equation}
\begin{equation}
\mbox{ and } \hspace*{0.1cm} P \big( |C_{s,x}|=\infty \big)=0. 
\end{equation}
\end{defi}

For example, consider a configuration where the size of all the clusters is smaller than an integer $m$, or a  configuration obtained by sub-critical or critical site percolation on $\zd$. Then, a $\ffp$ with one of these initial configurations satisfy a $\mbox{CCSB}(0,\delta,m)$ condition.

A natural question to ask is whether a $\ffp$ satisfies a $\mbox{CCSB}$ condition. If the forest-fire process has a CCSB at time $0$, then the following theorem states that after a certain time the answer is yes. The uniformity in time and in $\delta$ will be useful in the proof of the result.

\begin{theo}(Theorem 2 in \cite{Dur-these}) \label{theo ccsb}
Consider a $\ffp$ $\bar{\eta}$ on $G$ with parameter $\lambda>0$ that satisfies $CCSB(0,\frac{\lambda}{4d_{G}^2},m)$ for some integer $m$. Then for all real number $\gamma>0$, for all $\delta \in ]0,1]$, there exists an integer $ m_{\gamma,\lambda,d_G}(\delta)$ such that for all time $s \geq \gamma$, $\bar{\eta}$ has $CCSB(s,\delta,m_{\gamma,\lambda,d_G}(\delta))$.
\end{theo}

In \cite{Dur-these} an explicit formula for the parameter $m_{\gamma,\lambda,d_G}(\delta)$ appearing in the CCSB condition is given.

\vspace{0.3cm}

\subsubsection{Blur process}

\vspace{0.3cm}

We turn now to the blur-process. The aim of this process is to control the influence of long-distant sites on a given site. This process is associated with a forest-fire process $\bar{\eta}$, and defined from an initial time $t_0$ and a finite set $S$. The blur process indicates if the state of a site $x \in \bar{S}$ at time $t \geq t_0$ might be influenced by the state of the sites located outside the set $S$ at time $t_0$. If it might, the site is said to be $(t_0,S)$-blurred at time $t$. The word ``might" is important here. It means that the blur-process is a domination of the process of the real influence from outside $S$. If the site $x$ is not $(t_0,S)$-blurred at time $t$, its state at time $t$ can be recovered by using only the configuration at time $t_0$ and the Poisson processes of the sites in the set $S$.

The value of the $(t_0,S)$-blur process at a site $x \in \bar{S}$ at a time $t$ is either $2$ if the site is blurred or $0$ otherwise. Once a site is blurred, it is blurred forever. The set of blurred sites at $t_0$ is $N(S)$ and the clusters in $\eta(t_0)$ intersecting the boundary of $S$. As time evolves, this set grows toward the site $x$, due to births of trees that create new occupied paths. Here is the formal definition of the blur process.

\begin{defi}(blur process)
Let $t_0 \geq 0$ and $S \subset V$ a finite subset. The $(t_0,S)$-blur process $(\beta_{s,x})_{s \geq t_0, x \in \overline{S}}$ is a process with values in $\ \left\lbrace 0 , 2 \right\rbrace^{\overline{S}} $, right continuous that has the following properties. For all site $x$ in $\overline{S}$,
\begin{itemize}
\item the site $x$ is $(t_0,S)$-blurred \textbf{at time $t_0$} if and only if its cluster is connected to the boundary of $S$ : 
\[ \beta_{t_0,x}=\left \{ \begin{array}{l}  2 \mbox{  if } \hspace*{0.15cm} \overline{C}_{s,x} \cap N(S) \neq  \varnothing \\
0 \mbox{  else, }  \end{array} \right. \]

\item once a site is $(t_0,S)$-blurred, it remains $(t_0,S)$-blurred forever : 
\[ \mbox{  if } \beta_{s,x}=2, \mbox{  then } \forall s' \geq s,  \beta_{s',x}=2,\]

\item the site $x$ is $(t_0,S)$-blurred \textbf{at time $s$} $> t_0$ if and only if the set $\overline{C}_{s,x}$ contains a site that has been $(t_0,S)$-blurred before time $s$ : 
\[ \left\lbrace \beta_{s,x}=2 \right\rbrace = \left\lbrace \exists z \in \overline{C}_{s,x} \cap\overline{S} : \beta_{s^-,z}=2 \right\rbrace. \]

\end{itemize}  

\end{defi}

\vspace{0.3cm}

Like for the configuration process, the notation $\beta(t)$ will be used for $\lbrace \beta_{t,x} \rbrace_{x \in \zd}$.

\vspace{0.3cm}

It seems that the bigger the set $S$ is, the longer it takes for the site $x$ to be $S$-blurred. What happens when the size of the set $S$ tends to infinity? Is there a non empty interval during which a site is not influenced by ``infinity"? The answer is given by the following proposition, assuming a $\mbox{CCSB}$ condition is satisfied by the $\ffp$.

\begin{prop}(Proposition 2 in \cite{Dur-these}) \label{prop2 Durre}
For all $m \in \mathbb{N}$, there exists $\epsilon_m>0$ with the following property. Let $t \geq 0$ and suppose that the $\ffp$ has $\mbox{CCSB}(t,\frac{\lambda}{4d_G^{2}},m)$. Then for all sites $x$ in $V$, as $n$ tends to infinity, the probability that the site $x$ is $(t,B_n)$-blurred at time $t+\epsilon_m$ tends to zero : $$\lim_{n \rightarrow \infty} P( \beta_{t+\epsilon_m,x}=2)=0.$$ 
\end{prop}

\begin{re}
The real number $\epsilon_{m} >0$ is chosen such that $P(G_{0,\epsilon_{m}}>0)<\frac{1}{4 m d_G}$, where $G_{0,\epsilon_{m}}$ is the event that there is at least an occurrence in $[0,\epsilon_{m}]$ of a Poisson process with intensity $1$.
\end{re}

\begin{re}
This proposition is slightly different from the one written in \cite{Dur-these}. The additionnal hypothesis of almost-sure convergence at time $t$ is omitted here. Actually Dürre has shown that it is always satisfied for a $\ffp$ having a $CCSB$ condition at time zero.
\end{re}

\section{A candidate distribution}

To prove the existence of a stationary distribution, we exhibit a distribution and then show that it is stationary. The aim of this section is to explain how the candidate distribution is obtained. The idea is to use forest-fire processes on finite graphs. To get the translation-invariance property of the stationary distribution, a modified version of the boxes $B_k$ will be used. 

\subsection{Definition of ``finite volume" forest-fire processes}
\label{finitevolume} 

\vspace{0.3cm}

Let $k$ be a non negative integer. 

\vspace{0.3cm}

The first step is to define a graph. Contrary to what was mentionned above, the graph defined below is infinite. The fact that the process defined on it will behave like a finite-volume forest-fire explains why it is called ``finite volume''.   

\vspace{0.3cm}

Define a graph $G_k=(V^{(d)},E_k)$. The set of sites is $V^{(d)}=\zd$. The set of edges $E_k$ is bigger than the set of edges of the graph $\zd$. Some edges constructed from the box $B_k$ defined in Section~\ref{notations} will be added. The goal is to make opposite sites of the interior boundary of $B_k$ neighbors. 

Recall that $E^{(d)}$ is the set of edges corresponding to the hypercubic lattice of $\zd$ (see Section~\ref{notations}).

Define now the additionnal set of edges $A_k$. For all indices $j$ in $\{1, \dots , d\}$, consider a site with the $j$-th coordinate equal to $k$ and the other coordinates in the set $\lbrace -k, \dots , k \rbrace$. It is an interior boundary point. Then consider the opposite boundary site, i.e. the site with its $j$-th coordinate equal to $-k$ and the other coordinates exactly the same as the first site. Then $A_k$ consists in the edges between such pair of points for the whole interior boundary of $B_k$. More formally, for each $j \in \{1, \dots , d\}$ the edges between the following couples of sites are in $A_k$ : 

\[ \big \lbrace (i_1, \dots ,i_j,\dots, i_d) ; (i_1, \dots ,-i_j,\dots, i_d) \big \rbrace \mbox{ with }
\left \lbrace \begin{array}{l}
i_l=k \ \mbox{ if } \ l=j \\
i_l \in  \lbrace -k, \dots , k \rbrace \ \mbox{ if } \ l \neq j.
\end{array} \right. \]

Then the set $E_k$ is the union of the sets $E^{(d)}$ and $A_k$.

In Figure~\ref{graphGk}, a part of the additional edges for the graph $G_2$ restricted to the box $B_2 \subset \mathbb{Z}^2$ are drawn.

\begin{figure}[!h]
\centering
\includegraphics[scale=0.35]{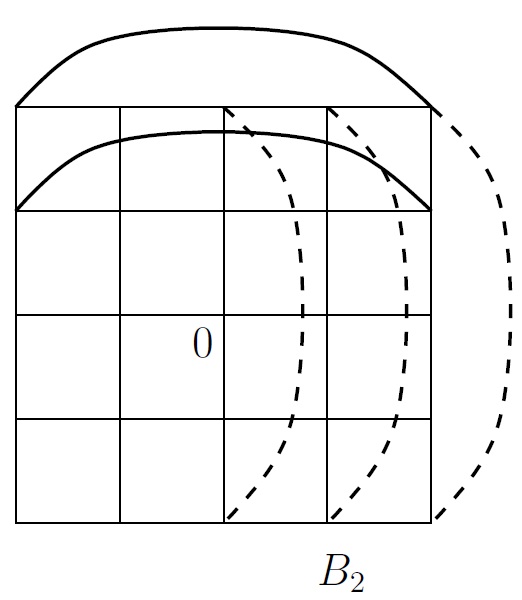}
\caption{Exemple of additional edges for the graph $G_2$}
\label{graphGk}
\end{figure}

Compared to the graph $\zd$, at most $d$ edges to each interior boundary site of $B_k$ is added in the graph $G_k$. So the vertex degree of $G_k$ is uniformly bounded by $d_f=3d$, instead of $2d$ for $\zd$. Let us denote $d_f=3d$ this bound.

The extra edges are not necessary to prove the existence of a stationnary distribution but there are used to prove that the stationnary distribution is translation-invariant.

\vspace{0.3cm}

Now the graph is defined, the second step is to consider a particular $\ffp$ on this graph.

\begin{defi}
A $G_k$-$\ffp$ $\bar{\eta}^{k}$ with parameter $\lambda>0$ is a $\ffp$ on the graph $G_k$ with parameter $\lambda>0$ with all the sites outside $B_{k}$ always empty : for all times $t \geq 0$, for all sites $x \in B_{k}^{c}$, $\eta^{k}_{t,x}=0$. 
\end{defi}

This forest-fire process behaves like a forest-fire process on $B_k$ where the opposite sites of the interior boundary are neighbors. The dynamic does not take into account the Poisson processes corresponding to the sites located outside $B_k$. However, its configurations take values in $\ \left\lbrace 0 , 1 \right\rbrace^{\zd} $. The $G_k$-$\ffp$ can be seen as a $\ffp$ on a discrete torus of dimension $d$.
  
\subsection{Construction of the candidate measure} 
\label{candidate} 

\vspace{0.3cm}

For each non negative integer $k$, let $\bar{\eta}^k$ be a $G_k$-$\ffp$. The configuration process $\eta^k$ is a finite state-space Markov chain. Consider two possible states of this Markov chain $\xi_1$ and $\xi_2$ and study the probability that the state $\xi_1$ leads to $\xi_2$. Consider the following event: all the occupied sites in $\xi_1$ are ignited and then a tree grows on each site which is occupied in $\xi_2$. As the growths and ignitions are driven by independent Poisson processes, this event occurs with a non negative probability. Therefore, this Markov chain is irreducible, recurrent and aperiodic. Thus it has a unique invariant distribution, which will be denoted by $\mu^k$. 

Since the space $\ \left\lbrace 0 , 1 \right\rbrace^{\zd} $ is compact, the sequence $(\mu^k)_{k \geq 0}$ has a weakly convergent subsequence. Let $\mu$ be its limit and $\mathcal{K}$ the set of indices of the subsequence. 

\begin{theo}
\label{theo mu}
Let $\eta$ be a $\ffp$ on $\zd$ with parameter $\lambda>0$. Then $\mu$ is an invariant and translation-invariant probability measure for $\eta$.
\end{theo}

As mentionned above, the $G_k$-$\ffp$ behaves exactly as a $\ffp$ on a discrete torus. Thus each distribution $\mu^k$ can be seen as a measure on a discrete torus, where we have the rotation-invariance property. Here a translation on the graph $G_k$ corresponds to a rotation on the discrete torus with $(2k+1)^d$ points, due to the $G_k$-$\ffp$ behavior. It follows that the weak limit $\mu$ is translation-invariant.

\section{Proof of the result}
\label{proof}

The goal of this section is to prove Theorem~\ref{main th}. It follows from Theorem~\ref{theo mu} stated above. This section is divided in three parts. The first part contains the proof of Theorem~\ref{theo mu} assuming two lemma. The lemma are then proved in the second and third part. 

\subsection{Proof of Theorem~\ref{theo mu}}
\label{proof sub}

To show that the distribution $\mu$ is stationary, the idea is to couple a $\ffp$ with initial configuration with law $\mu$ with $G_k$-$\ffps$ with initial configuration with law $\mu^k$ using the same Poisson processes. Then, the evolution of the processes will be compared using Dürre's tools.

A diagram of the main steps of the proof is given at the end of this first subsection.

Let $\lambda$ and $\gamma$ be two non negative real numbers. 

For simplification, as $\lambda$ and $d_f$ are fixed here, let us denote by $m_\gamma$ the real $m_{\gamma, \lambda, d_f}(\frac{\lambda}{4 \ d_f^2})$ of Theorem~\ref{theo ccsb} (recall that $d_f=3d$, see Section \ref{finitevolume}).

Let $\epsilon >0$ be a real number such that $P(G_{0,\epsilon}>0)<\frac{1}{4 m_\gamma d_f^2}$, where $G_{0,\epsilon}$ is the event that there is at least an occurrence in $[0,\epsilon]$ of a Poisson process with intensity $1$.

The aim is to prove the following : for all times $t < \epsilon $, for all cylinder events $A$, 
\begin{eqnarray}
\label{inv}
P^{\mu}(\eta(t) \in A)=\mu(A).
\end{eqnarray}

Let $t$ be a non negative real number such that $t < \epsilon $, and let $A$ be a cylinder event.

\vspace{0.3cm}

Throughout this proof, two specific sets of sites will be used. 

Since $A$ is a cylinder event, it is determined by the configuration on a finite number of sites. Let $I$ be the smallest box centered at the origin such that the event $A$ can be described using only sites in $I$. There exists $r_I$ such that $I=B_{r_I}$. 

We will need to work with a set containing $I$ but contained in $B_k$. Let $L$ be a non negative integer and $J$ be the box $B_{r_I+L}$. In order to simplify the notation, we will not write the dependence of $J$ on $L$.
%
%

\begin{figure}[H]
\centering
\includegraphics[scale=0.35]{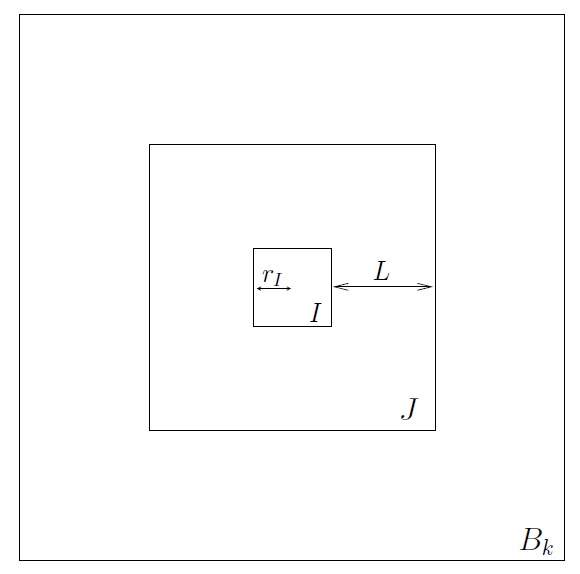}
\caption{Sets of sites}
\label{boxJ}
\end{figure}

For the remainder of the proof, only integers $k \in \mathcal{K}$ satisfying the relation $k>r_I+L$ will be considered.

\vspace{0.3cm}

Using the construction of the candidate measure $\mu$, a first upper bound is obtained : 

\begin{equation}
|P^{\mu}(\eta(t) \in A)- \mu(A)| \leq  |P^{\mu}(\eta(t) \in A)- \mu^k(A)| + |\mu^k(A)- \mu(A)|.
\end{equation}

Since $\mu^k$ is an invariant distribution for a $G_k$-$\ffp$, the equality (\ref{inv}) is satisfied for $\mu^k$. Turn now to study of the difference $|P^{\mu}(\eta(t) \in A)- P^{\mu^k}(\eta^k(t) \in A)|$.

\vspace{0.3cm}
 
Consider a $\ffp$ on $\zd$ and a $G_k$-$\ffp$, driven by the same Poisson processes on each site, and with initial configurations which coincide on $J$. Is it possible to find a time $t$ when the two configurations on a site $x$ in $I$ are not the same? The answer is yes if (and only if) the configuration outside the box $J$ influences the configuration in the box $I$ during the time interval $[0,t]$ (through some fires outside $J$). Which events can allow influence on $I$ from outside $J$? If at a time $t_0$ in $[0,t]$ there exists a path of occupied sites linking a site in $I$ to a site outside $J$, then a lightning falling in $J^c$ can burn trees in $I$. That means, using the blur process, the existence of a site in the set $I$ that is $(t,J)$-blurred.

\vspace{0.3cm}

However, in the difference we are interested in, the initial configuration of the two processes have the respective distribution $\mu$ and $\mu^k$. In order to use our previous argument, we are going to optimally couple the measures $\mu_{|J}$ and $\mu^{k}_{|J}$. The coupling of $\mu_{|J}$ and $\mu^{k}_{|J}$ is said to be optimal here if it maximizes the weight of the diagonal of the product of the spaces of configurations. For two configurations $\xi_1$ with law $\mu_{|J}$ and $\xi_2$ with law $\mu^{k}_{|J}$, the optimal coupling $\nu$ satisfies $\displaystyle{\nu(\xi_1 \neq \xi_2)=\inf_{\pi \in \mathcal{M}} \pi(\xi_1 \neq \xi_2)}$, where $\mathcal{M}$ is the set of all possible couplings of $\mu_{|J}$ and $\mu^{k}_{|J}$.

\vspace{0.3cm}
 
Now we couple a $\ffp$ $\bar{\eta}$ on $\zd$ with an initial configuration of law $\mu$ and the $G_k$-$\ffp$ $\bar{\eta}^k$ with an initial configuration of law $\mu^k$, by using the same Poisson processes and by optimally coupling $\mu_{|J}$ and $\mu^{k}_{|J}$. Then the following bound is obtained.

\begin{lemme}
\label{lemme majo-par-blur} 
Let us denote by $\beta$ the $(0,J)$-blur process associated with $\eta$. Then, 
\begin{eqnarray}
|P^{\mu}(\eta(t) \in A)- P^{\mu^k}(\eta^k(t) \in A)| \leq |I| \ \sup_{x \in I} P^{\mu}(\beta_{t,x}=2) + 2 \hspace*{0.1cm} d_{TV}(\mu_{|J},\mu^{k}_{|J}).
\end{eqnarray}
\end{lemme}

\vspace{0.3cm}

This lemma will be proved in the following section. The idea of the proof is to use the event ``the initial configurations of the two processes coincide on $J$". When it occurs, an upper bound is found using the blur process. Otherwise, the upper bound is given by the total variation between the two measures $\mu_{|J}$ and $\mu^{k}_{|J}$. 
 
\vspace{0.3cm}

Using this lemma, the upper bound of the desired quantity is : 
\begin{equation}
\label{borne fin}
|P^{\mu}(\eta(t) \in A)- \mu(A)| \leq |I| \ \sup_{x \in I} P^{\mu}(\beta_{t,x}=2) + 2 \hspace*{0.1cm} d_{TV}(\mu_{|J},\mu^{k}_{|J})+ |\mu^k(A)- \mu(A)|.
\end{equation}

\vspace{0.3cm}

The next step is to show that the upper bound of (\ref{borne fin}) tends to zero when $k$ and $L$ tend to infinity.

The distribution $\mu$ is the weak limit of a subsequence of $(\mu_k)_{k \geq 0}$. We first let $k$ go to infinity along the set of indices $\mathcal{K}$. Since $k> \mbox{radius}(J)=r_I+L$, $L$ does not depend on $k$ and $J$ is a finite set of sites, so : 
\begin{eqnarray}
&\displaystyle{\lim_{k \rightarrow \infty} |\mu^k(A)- \mu(A)|=0}& \\
&\displaystyle{\lim_{k \rightarrow \infty} d_{TV}(\mu_{|J},\mu^{k}_{|J})=0}.&
\end{eqnarray}

Recall that $\beta$ is a $(0,J)$-blur process. So the event $\lbrace \beta_{t,x}=2 \rbrace $ depends on $L$. To conclude we need to show that $\displaystyle{\hspace{0.1cm} \lim_{L \rightarrow \infty} P^{\mu}(\beta_{t,x}=2)=0}$. This is given by Proposition~\ref{prop2 Durre} applied with $t=0$ and $\delta=\frac{\lambda}{4d_f^{2}}$, provided the CCSB condition is satisfied. The purpose of the following lemma is to show that this hypothesis holds true.

\vspace{0.3cm}

\begin{lemme} \label{lemme mu}
For all non negative real number $\delta$, a $\ffp$ on $\zd$ with parameter $\lambda>0$, with an initial configuration with distribution $\mu$, has $\mbox{CCSB}(0,\delta,m_{\gamma,\lambda,d_f}(\delta))$.
\end{lemme}

\begin{re}
In the last parameter of the $\mbox{CCSB}$ condition, the bound used for the vertex degree of $\zd$ is $d_f=3d$ and not $2d$ which is the optimal bound. This is due to our argument to prove this lemma, which uses the $G_k$-$\ffp$. Nevertheless this is sufficient here to show our theorem.
\end{re}

Using Lemma~\ref{lemme mu} with $\delta=\frac{\lambda}{4d_f^{2}}$, we claim that a $\ffp$ with an initial configuration with law $\mu$ has $\mbox{CCSB}(0,\frac{\lambda}{4d_f^{2}},m_{\gamma,\lambda,d_f}(\frac{\lambda}{4d_f^{2}}))$. 

\vspace{0.3cm}

By Proposition~\ref{prop2 Durre}, there exists $\epsilon_{m_\gamma}$ such that for all site $x$ in $\zd$, as $L$ tends to infinity, the probability that the site $x$ is $(0,J)$-blurred at time $\epsilon_{m_\gamma}$ tends to zero. That is to say, with the choice of $\epsilon$ made here, that for all $t<\epsilon$, for all $x \in \zd$ $\displaystyle{\hspace{0.2cm} \lim_{L \rightarrow \infty} P^{\mu}(\beta_{t,x}=2)=0}$. 

Since $|I| < \infty$, this concludes the proof. \qed

\vspace{0.2cm}

We conclude this part with a diagram of the main steps of the proof.

\vspace{0.1cm}

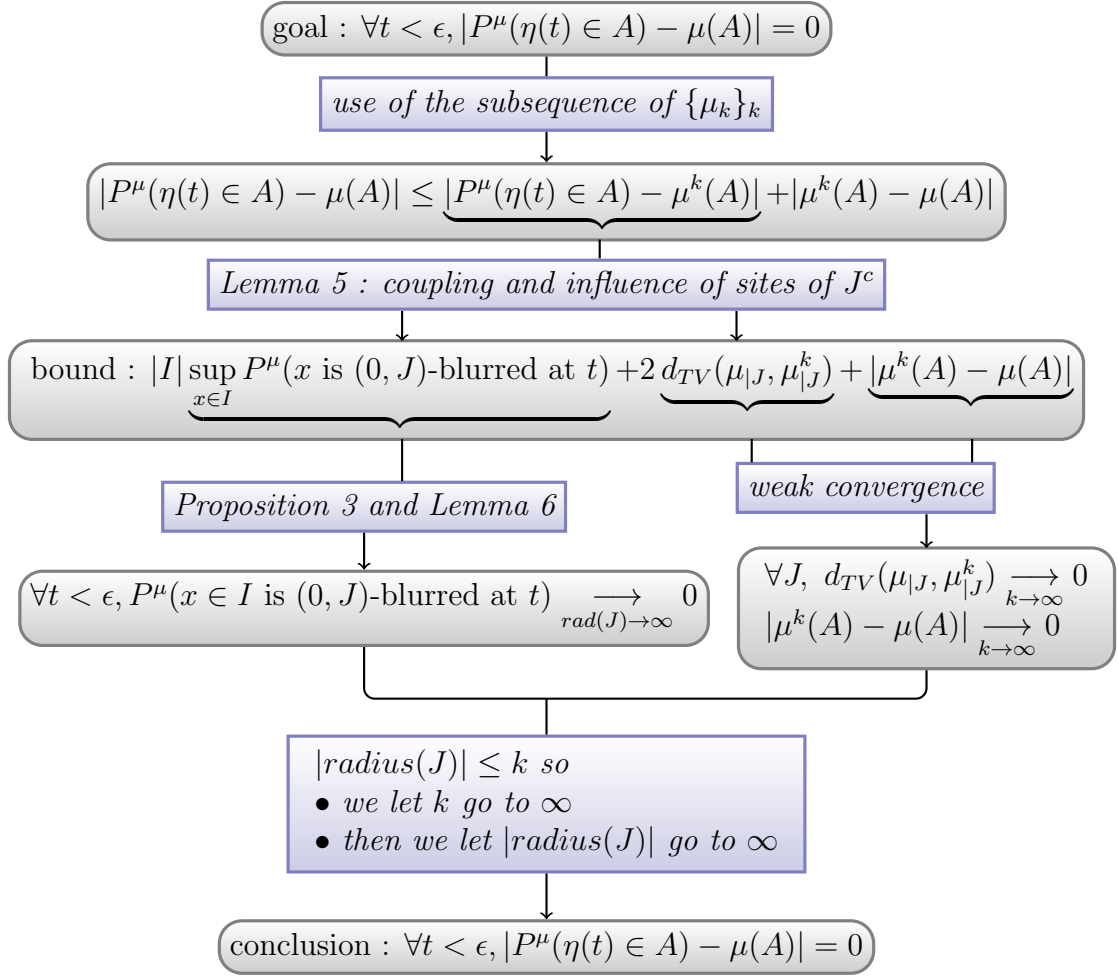
\begin{figure}[H]
\centering

\begin{tikzpicture}[start chain,node distance=5mm,every join/.style={->} ]
\node (rond1) at (0,7.2) [terminal] {goal : $\forall t<\epsilon, |P^{\mu}(\eta(t) \in A) - \mu(A)| = 0$  };
\node (rect1) at (0,6.2) [nonterminal] {use of the subsequence of $\lbrace \mu_k \rbrace_k$};
\node (rond2) at (0,4.9) [terminal] {$|P^{\mu}(\eta(t) \in A)- \mu(A)| \leq  \underbrace{|P^{\mu}(\eta(t) \in A) - \mu^k(A)|} + |\mu^k(A)- \mu(A)| $};
\node (rect2) at (0,3.8) [nonterminal]{
Lemma~\ref{lemme majo-par-blur} : coupling and influence of sites of $J^c$};
\node (rond3) at (0,2.4) [terminal] {$\displaystyle{\mbox{ bound : } |I| \underbrace{\sup_{x \in I} P^\mu( x \mbox{ is } (0,J)\mbox{-blurred at } t)} + 2 \underbrace{d_{TV}(\mu_{|J},\mu^{k}_{|J})}  + \underbrace{|\mu^k(A)- \mu(A)| } }$};
\node (expl1) at (-2.4,0.85) [nonterminal]{
Proposition~\ref{prop2 Durre} and Lemma~\ref{lemme mu}};
\node (res1) at (-2.4,-0.50) [terminal] {$\forall t<\epsilon, P^\mu( x \in I \mbox{ is } (0,J)\mbox{-blurred at } t) \underset{rad(J) \to \infty}{\longrightarrow} 0$};
\node (expl2) at (4.2,1.1) [nonterminal]{weak convergence};
\node (res2) at (5,-0.5) [terminal] {\begin{tabular}{l}
$\forall J, \ d_{TV}(\mu_{|J},\mu^{k}_{|J}) \underset{k \to \infty}\longrightarrow 0$ \\
$|\mu^k(A)- \mu(A)|\underset{k \to \infty}\longrightarrow 0 $
\end{tabular}};
\node (rect3) at (0,-3.1) [nonterminal] {\begin{tabular}{l}
$|radius(J)| \leq k$ so \\
$\bullet$ we let $k$ go to $\infty$\\
$\bullet$ then we let $|radius(J)|$ go to $\infty$
\end{tabular}};
\node (rond4) at (0,-5) [terminal] {conclusion : $\forall t<\epsilon, |P^{\mu}(\eta(t) \in A) - \mu(A)| = 0$  };
\path (rond1) edge[thick,-] (rect1) ; 
\path (rect1) edge[thick,->] (rond2) ; 
\path (0.7,4.4) edge[thick,-] (0.7,4.15); 
\path (2.5,3.05) edge[thick,<-] (2.5,3.45);
\path (-1.9,3.05) edge[thick,<-] (-1.9,3.45);
\path (-1.9,1.19) edge[thick,-] (-1.9,1.75);
\path (2.7,1.44) edge[thick,-] (2.7,1.75);
\path (5.6,1.44) edge[thick,-] (5.6,1.75);
\path (rect3) edge[thick,->] (rond4);
\path (expl1) edge[thick,->] (res1);
\path (5,0.3) edge[thick,<-] (5,0.75);
\draw [thick,-,rounded corners] (node cs:name=res1,anchor=south) |- (0,-1.7);
\draw [thick,-,rounded corners] (node cs:name=res2,anchor=south) |- (0,-1.7);
\draw [thick] (node cs:name=rect3,anchor=north) |- (0,-1.7);

\end{tikzpicture}

\caption{Scheme of the proof}
\label{globalscheme}
\end{figure}


\subsection{Proof of Lemma~\ref{lemme majo-par-blur}}

\label{remaining}  
 
As mentionned in the proof of Theorem~\ref{theo mu}, the idea to prove this lemma is to use the event ``the initial configurations of the two processes coincide on $J$". Then, the study differs depending on whether this event occurs or not. A diagram of the main steps of the proof is given at the end of the section. 
 
\vspace*{0.3cm}

Recall that we coupled a $\ffp$ $\eta$ on $\zd$ with an initial configuration of law $\mu$ and the $G_k$-$\ffp$ $\eta^k$ with an initial configuration of law $\mu^k$, by using the same Poisson processes on each site and by optimally coupling $\mu_{|J}$ and $\mu^{k}_{|J}$. Let $\nu$ be the distribution of the optimal coupling. We have 
\[
|P^{\mu}(\eta(t) \in A)- P^{\mu^k}(\eta^k(t) \in A)| = |P^{\nu}((\eta(t),\eta^k(t)) \in A \times X)-P^{\nu}((\eta(t),\eta^k(t)) \in X \times A)|.
\]
Consider now the event ``the initial configurations of the two processes coincide on $J$". Let $C_J=\lbrace \eta(0)_{|J} = \eta^k(0)_{|J} \rbrace$ be this event. Denote by $C_{J}^c$ the complement of this event. Then we can upper bound the difference we are interested in : 

\begin{equation}
|P^{\mu}(\eta(t) \in A)- P^{\mu^k}(\eta^k(t) \in A)|  =  |Dif+Eq|  \leq  |Dif|+|Eq| 
\end{equation} 

with
\begin{eqnarray}
Dif & = & P^{\nu} \big(  \left[  (\eta(t),\eta^k(t)) \in A \times X \right]  \cap C_{J}^c \big) 
 -  P^{\nu}\big(  \left[  (\eta(t),\eta^k(t)) \in X \times A \right]  \cap  C_{J}^c \big) \hspace{0.1cm}
\end{eqnarray}
\begin{eqnarray}
\label{equa eq}
Eq  & = & P^{\nu}\big(  \left[  (\eta(t),\eta^k(t)) \in A \times X \right]  \cap C_J \big) 
 - \hspace*{0.2cm} P^{\nu}\big(  \left[  (\eta(t),\eta^k(t)) \in X \times A \right]  \cap C_J \big).\hspace{0.1cm}
\end{eqnarray}

\vspace{0.3cm}

$\star$ Study of $Dif$ : 

Since $\mu_{|J}$ and $\mu^{k}_{|J}$ are optimally coupled,
\[P^{\nu}(\eta(0)_{|J} \neq \eta^k(0)_{|J}) \leq  \hspace*{0.1cm} d_{TV}(\mu_{|J},\mu^{k}_{|J}) \]

\noindent and we get 

\begin{equation}
\label{dif}
|Dif| \leq 2 \hspace*{0.1cm} P^\nu(\eta(0)_{|J} \neq \eta^k(0)_{|J}) \leq 2 \hspace*{0.1cm} d_{TV}(\mu_{|J},\mu^{k}_{|J}).
\end{equation}

\vspace{0.3cm}

$\star$ Study of $Eq$ : 

After using the Bayes formula in (\ref{equa eq}), we are interested in the difference 
\begin{equation}
Eq = \big | P^{\nu}\big( \ (\eta(t),\eta^k(t)) \in A \times X  \ \big | \ C_J \big) - P^{\nu}\big( \ (\eta(t),\eta^k(t)) \in X \times A  \ \big | \ C_J \big) \big | \cdot P(C_J).
\end{equation}

\noindent We can use indicator functions to rewrite the difference in the righthand term : 
\[  Eq = \big | E^{\nu}(1_{X \times A}(\eta(t),\eta^k(t)) -1_{A \times X}(\eta(t),\eta^k(t)) \ | \ C_J ) \big | \cdot P(C_J) . \]

When is the difference $1_{X \times A}(\eta(t),\eta^k(t)) -1_{A \times X}(\eta(t),\eta^k(t))$ non zero? It occurs only when $(\eta(t),\eta^k(t)) \in {A^c \times A}$ or when $(\eta(t),\eta^k(t)) \in {A \times A^c}$. These events occur when there exists a time $t_1>0$ and a site $x$ in $I$ where the two configurations are not the same, i.e. $\eta_{t_1,x} \neq \eta^{k}_{t_1,x}$. The two $\ffps$ are driven by the same Poisson processes and we conditionned by the event ``their initial configurations coincide on $J$". So this non equality must be the consequence of the influence of at least one cluster outside $J$ on the box $I$. This is equivalent to the existence of a site in $I$ whose state is influenced by the state of a site outside $J$. More precisely, there exists a site $x$ in $I$ which is $(0,J)$-blurred at time $t_1$. 

\vspace{0.3cm}

Firstly we give some notations associated with the blur processes. We denote by $\beta$ the $(0,J)$-blur process associated with $\bar{\eta}$ and by $\beta^k$ the $(0,J)$-blur process associated with $\bar{\eta}^k$. To simplify the notation, we do not write the dependence in $J$ (neither in $L$) of the blur processes. Let $NB$ be the set of configurations where all the sites of the box $I$ have the value $0$ and $NB_x$ be the set of configurations where the site $x$ has the value $0$.

\vspace{0.3cm}

Then with these notations we can write
 \[|1_{X \times A}(\eta(t),\eta^k(t)) - 1_{A \times X}(\eta(t),\eta^k(t))| \leq 1_{NB^c  \times NB^c } (\beta(t),\beta^k(t)) .\] So
\begin{eqnarray*}
|Eq| & \leq & E^{\nu}(1_{NB^c  \times NB^c } (\beta(t),\beta^k(t)) \ | \ C_J ) \hspace*{0.2cm} P^{\nu}(C_J) \\
& \leq & P^\nu( \ (\beta(t),\beta^k(t)) \in NB^c  \times NB^c \ | \ C_J) \cdot P^{\nu}(C_J) \\
\end{eqnarray*}

The $(0,J)$-blur process associated with a $\ffp$ is defined only with the Poisson processes on the sites of the set $J$ and the initial configuration of the sites of the box $J$. We know that the two initial configurations coincide on $J$ and that the two processes are driven by the same Poisson processes on each site. Thus, the blur processes $\beta$ and $\beta^k$ are the same. So \[|Eq| \leq P^\mu( \beta(t) \in NB^c).\] 

\vspace{0.3cm}

The complements of the events $NB$ and $NB_x$ satisfy the relation $\displaystyle{ NB^c= \bigcup_{x \in I} NB^{c}_{x}},$ so 
\begin{equation}
\label{eq}
|Eq| \leq |I| \hspace{0.1cm} \sup_{x \in I} P^{\mu}(\beta(t) \in NB^{c}_{x}).
\end{equation}

$\star$ Conclusion 

\vspace{0.3cm}

Since $P^{\mu}(\beta(t) \in NB^{c}_{x})=P^{\mu}(\beta_{t,x}=2)$, using equations (\ref{dif}) and (\ref{eq}) we obtain the desired result.

\[ |P^{\mu}(\eta(t) \in A)- P^{\mu^k}(\eta^k(t) \in A)| \leq 2 \hspace*{0.1cm} d_{TV}(\mu_{|J},\mu^{k}_{|J}) + |I| \hspace{0.1cm} \sup_{x \in I} P^{\mu}(\beta_{t,x}=2). \]
\qed

\begin{figure}[H]

\begin{tikzpicture}[start chain,node distance=5mm,every join/.style={->} ]
\node (res1) at (-1,5.5) [terminal] {\begin{tabular}{c} goal : \\ $\displaystyle{|P^{\mu}(\eta(t) \in A) - P^{\mu^k}(\eta^k(t) \in A)| \leq |I| \ \sup_{x \in I} P^\mu( x \mbox{ is } (0,J)\mbox{-blurred at } t) + 2 \hspace*{0.1cm} d_{TV}(\mu_{|J},\mu^{k}_{|J})}$\end{tabular} };
\node (justif2) at (-1,3.5) [nonterminal] {\begin{tabular}{l}
optimal coupling of $\mu$ and $\mu_k$  \\
$\rightsquigarrow$ either the 2 initial processes coincide on $J$ (Eq) \\
$\rightsquigarrow$ or they are different (Dif)
\end{tabular}};
\node (res3) at (-1,1.5) [terminal] {$|P^{\mu}(\eta(t) \in A)- P^{\mu^k}(\eta^k(t) \in A)| \leq |Eq|+|Dif|$};
\node (justEq1) at (-4,-0.2) [nonterminal] {study of the influence of sites outside $J$};
\node (resEq1) at (-4,-1.7) [terminal] {\begin{tabular}{lll}
$|Eq|$ & $\leq$ & $P^\mu(\exists x \in I : x \mbox{ is } (0,J)\mbox{-blurred at } t)$ \\
& $\leq$ & $|I| \ \sup_{x \in I} P^\mu( x \in I \mbox{ is } (0,J)\mbox{-blurred at } t)$
\end{tabular}};
\node (justDif1) at (4,-0.2) [nonterminal] {coupling property};
\node (resDif1) at (4,-1.7) [terminal] {$|Dif| \leq 2 \hspace*{0.1cm} d_{TV}(\mu_{|J},\mu^{k}_{|J})$};
\node (ccl) at (-1,-3.7) [terminal] {conclusion : the desired upper bound is obtained};
\path (res1) edge[thick,-] (justif2) ; 
\path (justif2) edge[thick,->] (res3) ; 
\draw [thick,-,rounded corners] (node cs:name=justEq1,anchor=north) |- (0,0.7);
\draw [thick,-,rounded corners] (node cs:name=justDif1,anchor=north) |- (0,0.7);
\draw [thick] (node cs:name=res3,anchor=south) |- (0,0.7);
\path (justEq1) edge[thick,->] (resEq1) ;
\path (justDif1) edge[thick,->] (resDif1) ;
\draw [thick,-,rounded corners] (node cs:name=resEq1,anchor=south) |- (0,-2.8);
\draw [thick,-,rounded corners] (node cs:name=resDif1,anchor=south) |- (0,-2.8);
\draw [thick,<-] (node cs:name=ccl,anchor=north) |- (0,-2.8);
\end{tikzpicture} 
 
\caption{Scheme of the proof of Lemma~\ref{lemme majo-par-blur}}
\label{secondstep}
\end{figure}
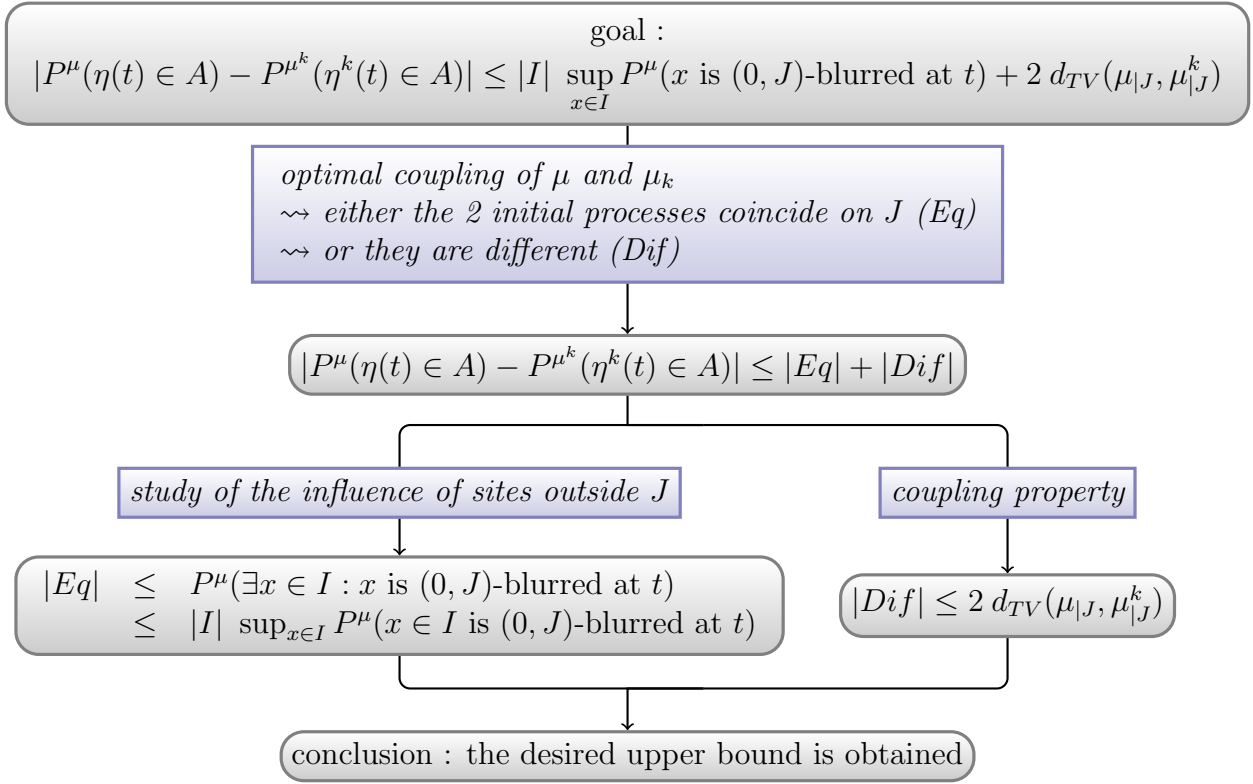

\subsection{Proof of Lemma~\ref{lemme mu}}

\label{studyblur} 

The construction of the distribution $\mu$ as a weak limit of the sequence $\lbrace \mu^k \rbrace_{k \in \mathcal{K}}$ is used to prove this lemma. We will need the CCSB condition for a $G_k$-$\ffp$, which is the result of the following lemma. 

\begin{lemme}  \label{lemme mu n}
Let $n$ be in $\mathbb{N}$ and  $\delta$ be a non negative real number. Consider a $G_n$-forest-fire process $\alpha^n$ with parameter $\lambda$, with an initial configuration with distribution $\mu^n$. 

Then $\alpha^n$ has $\mbox{CCSB}(0,\delta,m_{\gamma,\lambda,d_f}(\delta))$.
\end{lemme}

\begin{dem}(of Lemma~\ref{lemme mu n})

Let $n$ and $m$ be in $\mathbb{N}$.

\vspace{0.3cm}

The goal is to show a CCSB condition. Therefore let $B$ and $D$ be two subsets of $B_n$, $x$ a site in $B_n \setminus D$ and consider 
\[ \mbox{Cond}_{B,D}= \Big \lbrace \bigcup_{y \in B} C_{y}=D \Big \rbrace  \hspace*{0.5cm} \mbox{and} \hspace*{0.5cm}  \mbox{Size}_{x,m}=\big \lbrace|C_{x}|>m \big \rbrace.  \]

To simplify the notations, we denote $m_{\gamma,\lambda,d_f}(\delta)$ by $m_{\gamma,\delta}$ in this section.

\vspace{0.3cm}

The $G_n$-$\ffp$ $\alpha^n$ has $\mbox{CCSB}(0,\delta,m_{\gamma,\lambda,d_f}(\delta))$ if 
\begin{equation}
\label{prev}
P^{\mu^n}(\alpha^n(0) \in \mbox{Size}_{x,m_{\gamma,\delta}} \cap \alpha^n(0) \in \mbox{Cond}_{B,D} ) \leq \delta  P^{\mu^n}( \alpha^n(0) \in \mbox{Cond}_{B,D}). 
\end{equation}

Since $\mu^n$ is the stationary distribution of a $G_n$-$\ffp$, (\ref{prev}) can be written : 
\begin{equation}
\label{eq mu n}
\mu^n(\mbox{Size}_{x,m_{\gamma,\delta}} \cap \mbox{Cond}_{B,D} ) \leq \delta  \mu^n( \mbox{Cond}_{B,D}). 
\end{equation}

A forest-fire process having a CCSB condition at time $0$ satisfies one for all times $t>\gamma$ by Theorem~\ref{theo ccsb}. Thus to show (\ref{eq mu n}), we will use this result combined with the ergodic theorem.

\vspace{0.3cm}

Consider a $G_n$-$\ffp$ $\beta^n$ with parameter $\lambda>0$, with an initial configuration with law $\rho^n$ and having $\mbox{CCSB}(0,\frac{\lambda}{4d_f^{2}},m)$.

\vspace{0.3cm}

Since $\mu^n$ is the unique stationary distribution of a $G_n$-$\ffp$, the ergodic theorem gives
\begin{equation}
\label{th ergod}
 \mu^n(\mbox{Size}_{x,m_{\gamma,\delta}} \cap \mbox{Cond}_{B,D})=\lim_{t \rightarrow \infty} \frac{1}{t} \int_{0}^{t} P^{\rho^n}(\beta^{n}_{s} \in \mbox{Size}_{x,m_{\gamma,\delta}} \cap \beta^n(s) \in \mbox{Cond}_{B,D}) ds . 
 \end{equation}
  
Recall that $\gamma$ was set in Section~\ref{proof sub} and take $t>\gamma$.

By Theorem~\ref{theo ccsb}, for all times $s$ larger than $\gamma$, the $\ffp$ $\beta^n$ has $\mbox{CCSB}(s,\delta,m_{\gamma,\delta})$, so

\begin{equation*}
\forall s \geq \lambda, \ P^{\rho^n}(\beta^{n}_{s} \in \mbox{Size}_{x,m_{\gamma,\delta}} \cap \beta^n(s) \in \mbox{Cond}_{B,D}) \leq \delta  P^{\rho^n}(\beta^{n}_{s} \in \mbox{Cond}_{B,D}).
 \end{equation*}
 
Therefore, 
\begin{eqnarray*}
 \frac{1}{t} \int_{\gamma}^{t} P^{\rho^n}(\beta^{n}_{s} \in \mbox{Size}_{x,m_{\gamma,\delta}} \cap \beta^n(s) \in \mbox{Cond}_{B,D}) ds & \leq  & \frac{1}{t} \int_{\gamma}^{t} \delta \ P^{\rho^n}(\beta^{n}_{s} \in \mbox{Cond}_{B,D}) ds \\
 &  \leq  & \delta \ \frac{1}{t} \int_{0}^{t}  P^{\rho^n}(\beta^{n}_{s} \in \mbox{Cond}_{B,D}) ds.
\end{eqnarray*}

Since \[ \hspace*{0.2cm} \frac{1}{t} \int_{0}^{\gamma} P^{\rho^n}(\beta^{n}_{s} \in \mbox{Size}_{x,m_{\gamma,\delta}} \cap \beta^n(s) \in \mbox{Cond}_{B,D}) ds \leq \frac{\gamma}{t}, \] 

we get a bound for the integral :

\begin{equation*} 
\frac{1}{t} \int_{0}^{t} P^{\rho^n}(\beta^{n}_{s} \in \mbox{Size}_{x,m_{\gamma,\delta}} \cap \beta^n(s) \in \mbox{Cond}_{B,D}) ds \leq \frac{\gamma}{t} + \delta \frac{1}{t} \int_{0}^{t} P^{\rho^n}(\beta^{n}_{s} \in \mbox{Cond}_{B,D}) ds. 
\end{equation*}
 
Finally, we let $t$ tend to infinity and use (\ref{th ergod}) twice to get (\ref{eq mu n}). This concludes the proof. \qed
\end{dem}

We can now show the lemma concerning the distribution $\mu$.

\begin{dem} (of Lemma~\ref{lemme mu})

Let $\delta>0$.

The idea of the proof is to use the fact that the measure $\mu$ is a weak limit of a subsequence of $(\mu^k)_{k \geq 0}$, combined with the previous lemma.

To show the CCSB condition, consider two finite subsets $B$ and $D$ of $\zd$, a site $x$ in $\zd \setminus D$  and
\[ \mbox{Cond}_{B,D}= \Big \lbrace \bigcup_{y \in B} C_{y}=D \Big \rbrace  \hspace*{0.5cm} \mbox{and} \hspace*{0.5cm}  \mbox{Size}_{x,m_{\gamma,\delta}}=\big \lbrace|C_{x}|>m_{\gamma,\delta} \big \rbrace.  \]

For all finite subsets $B$ and $D$ of $\zd$ and all sites $x$ in $\zd \setminus D$, the sets $\mbox{Cond}_{B,D}$ and $\mbox{Size}_{x,m_{\gamma,\delta}}$ are defined by a finite number of sites. Therefore we can write 
\begin{eqnarray}
\label{lim}
\mu(\mbox{Size}_{x,m_{\gamma,\delta}} \cap \mbox{Cond}_{B,D}) & = & \lim_{k \rightarrow \infty} \mu^k(\mbox{Size}_{x,m_{\gamma,\delta}} \cap \mbox{Cond}_{B,D}) \\
\label{lim2}
 \mu(\mbox{Cond}_{B,D}) & = & \lim_{k \rightarrow \infty} \mu^k(\mbox{Cond}_{B,D}).
\end{eqnarray}

By the previous lemma 
\begin{equation} 
\mu^k(\mbox{Size}_{x,m_{\gamma,\delta}} \cap \mbox{Cond}_{B,D}) \leq \ \delta \ \mu^k(\mbox{Cond}_{B,D}).
\end{equation}

The result is obtained by letting $k$ tend to infinity :
\begin{equation*} 
\mu(\mbox{Size}_{x,m_{\gamma,\delta}} \cap \mbox{Cond}_{B,D}) \leq \ \delta \ \mu(\mbox{Cond}_{B,D}).
\end{equation*}
 \qed

\end{dem}


\vspace{0.5cm}

\sc{Alice Stahl}

\sc{Institut de Mathématiques de Toulouse}

\sc{Université de Toulouse, F-31062 Toulouse, France}

alice.stahl@math.univ-toulouse.fr

\end{document}